\documentclass[conference,letterpaper]{IEEEtran}
\usepackage[letterpaper, left=0.7in, right=0.7in, bottom=0.637in, top=0.833in]{geometry}
\IEEEoverridecommandlockouts
\usepackage[sorting=none,style=ieee]{biblatex}
\bibliography{LCSS_ref}

\usepackage{epstopdf}
\usepackage{amsmath,amssymb,amsfonts}

\usepackage{float}
\usepackage{xcolor}
\usepackage{algorithmic}
\usepackage{mathtools}
\usepackage{graphicx}
\usepackage{textcomp}
\usepackage{xcolor}
\usepackage{tabularx,booktabs}
\usepackage{varioref}
\usepackage{Keywords}
\usepackage{tikz-cd}
\usepackage{soul}
\newcolumntype{C}{>{\centering\arraybackslash}X} % centered version of "X" type
\def\BibTeX{{\rm B\kern-.05em{\sc i\kern-.025em b}\kern-.08em
    T\kern-.1667em\lower.7ex\hbox{E}\kern-.125emX}}
\newenvironment{@abssec}[1]{%
     \if@twocolumn
       \section*{#1}%
     \else
       \vspace{.05in}\footnotesize
       \parindent .2in
         {\upshape\bfseries #1. }\ignorespaces 
     \fi}
     {\if@twocolumn\else\par\vspace{.1in}\fi}

\begin{document}

\title{Uncovering Quasi-periodic Nature of Physical Systems: A Case Study of Signalized Intersections}

\author{Suddhasattwa~Das\textsuperscript{1},
        Shakib~Mustavee\textsuperscript{2},
        ~Shaurya~Agarwal\textsuperscript{2}
\thanks{$^{1}$The author is with the Department
of Mathematics, George Masion University, Fairfax, Virginia, USA  (email: iamsuddhasattwa@gmail.com)}% <-this % stops a space
\thanks{$^{2}$The authors are with the Department of Civil Engineering, University of Central Florida, Orlando, USA (email: smustavee@gmail.com), (email: Shaurya.Agarwal@ucf.edu)}}% <-this % stops a space
%\thanks{Manuscript received April 19, 2005; revised August 26, 2015.}}

\maketitle

%Main body starts

\begin{abstract}
This paper presents a novel approach to analyze quasiperiodically driven dynamical systems. It aims to develop a complete data-driven framework for modeling such unknown dynamics. To achieve this, we characterize Koopman eigenfrequencies as generating frequencies of the quasiperiodic driver of the system. We compute true eigenfrequencies of Koopman operators by applying the theory of Reproducing Kernel Hibert Space (RKHS) and results from ergodic theory. We also demonstrate the decomposition of quasiperiodically driven dynamics into two components, i) the quasiperiodic driving source with generating frequencies and ii) the driven nonlinear dynamics. A unique aspect of the proposed framework is that it applies to the analysis of systems where the periodic component is either non-dominant or even absent. As a case study, we analyze a system of nine traffic signalized intersections. The proposed framework accurately reconstructs the measured queue lengths of the signalized intersections and makes stable long-term predictions. 

\end{abstract}

\begin {IEEEkeywords}
 Koopman operator, Reproducing Kernel Hibert Space, Quasiperiodic systems, Signalized intersections
\end{IEEEkeywords}

\section{Introduction} \label{intro}
Many physical phenomena show periodicities, such as oscillations and limit cycles, characterized by a single frequency or period. There are also many phenomena, such as the arrival of fall foliage, the arctic ice cycle, and planetary dynamics, in which there are multiple driving frequencies. As a result of which the states of the system do not show exact and regular periodicity. Such systems are called \emph{quasiperiodic systems}. More complicated systems exhibit chaos and a high degree of nonlinearity but having a quasiperiodic component as a driving source. Examples are astronomical systems \cite{DasJim2017_SuperC}, climate systems \cite{VautardGhil89, SlawinskaGiannakis16}, geophysical flows on periodic domains \cite{GiannakisDas_tracers_2019}, epidemics \cite{mustavee2021linear}, and computational neuroscience \cite{marrouch2020data}. %We call this \emph {quasiperiodically driven} dynamics. 

The goal of this paper is to develop a robust method for reconstructing quasiperiodically driven dynamical systems using the Koopman operator framework.  Koopman operator is originally a linear operator that describes the spatiotemporal evolution of a nonlinear dynamical system. Projection of observables on the eigenfunction of the Koopman operator is defined as Koopman modes. Koopman modes can efficiently capture quasiperiodic components of a flow and outperform proper orthogonal decomposition characterizing the evolution on limit cycles, and tori \cite{arbabi2017study}. Koopman operator has also proven to be successful for describing highly non-periodic dynamics by decomposing the underlying dynamics into periodic and quasiperiodic patterns  \cite{marrouch2020data}.  In Section~\ref{sec:eigen}, we describe quasiperiodicity in terms of eigenfunctions of the Koopman operator. However, an accurate estimation of the eigenfunctions and eigenfrequencies has remained an elusive task. Techniques such as DMD \cite{arbabi2017study, shabab2021exploring}, deep neural networks \cite{YeungEtAl2019, GononOrtega_2019_a}, or Fourier averaging \cite{DasJim2017_SuperC} are inadequate when the dynamics has a substantial chaotic component. %For such systems, both DMD-type and Fourier techniques fail to identify the true eigenfrequencies. 

Our novel approach uses a technique from \emph{RKHS interpolation theory}  developed in \cite{DasGiannakis_RKHS_2018} to extract the true eigenfrequencies. After that, we proceed to identify different components of the dynamics in a manner similar in structure to \cite{marrouch2020data}. This approach leads to a robust theoretical and numerical method to reconstruct quasiperiodically driven systems having strong chaotic components by identifying their true eigenfrequencies. As a case study, we study the queue length dynamics on a corridor of traffic signalized intersections. Such a system has quasiperiodic driving sources, but the measurement data also reflects unpredictable phenomena such as accidents, road constructions, and human factors. It makes the study more challenging.

\textit{Contributions:} The contributions of this paper is as follows:
\begin{itemize}
\item We provide an alternative mathematical formulation of the quasiperiodic coordinates in terms of the \emph{Koopman operator} that correctly recovers both the periodic and chaotic part of the dynamics.% based on rigorous results from ergodic theory \cite{DasGiannakis_RKHS_2018}. 
\item The proposed framework is even applicable to the systems where the periodic component is either non-dominant or even absent.
\item We describe the physical significance of the quasiperiodic sources obtained from the proposed technique
\end{itemize}

\textit{Outline:} Rest of the paper is arranged as follows: The dynamical systems theory related to quasiperiodically driven systems is described in Section~\ref{dynamical_theory}. We discuss the relevant concepts and techniques in Section~\ref{sec:eigen} and Section~\ref{sec:kernel}. The proposed data-driven implementation is described in Section~\ref{sec:data_driven_algorithms}. Section~\ref{sec:Alfaya} performs a case study using real measurements from traffic intersections and discusses the results.

\section{Proposed Framework} \label{section1}

\subsection{Formulation}\label{dynamical_theory}

A quasiperiodically driven system has the form
\begin{equation} \begin{split} \label{eqn:def:quasi_driven}
\theta_{n+1} &= \theta_n + \vec\rho \bmod \mathbb{T}^d \\
x_{n+1} &= g( x_n, \theta_n)
\end{split},\end{equation}
%
%\red{$\theta \bmod 2\pi$ is the unique number $\hat{\theta}\in [0,2\pi)$ such that $\theta= \hat{\theta} + k 2\pi$ for some integer $k$.}
%
where $f$ is some nonlinear function. The angular coordinate $\theta$ is a point on a $d$-dimensional torus $\mathbb{T}^d$. $\theta$ represents the \emph{phase} of a driving quasiperiodic system. The vector $\vec\rho$ is called the \emph{rotation vector} \cite{Herman1, Arnold1965}, it represents the angular increments at each step for each of the coordinates of $\theta$. If the underlying system arises from a continuous time system by taking samples at intervals $\Delta t$, then $\vec\rho = \Delta t \vec\omega$ for some angular frequency vector $\vec\omega$. The variable $x$ lies in some abstract or unknown manifold $\mathcal{X}$. Let $\Omega := \mathbb{T}^d \times \mathcal{X}$. Thus \eqref{eqn:def:quasi_driven} is a \emph{one-way coupled} or \emph{skew-product} dynamical system on the space $\Omega$ :
\begin{equation} \label{eqn:def:quasi_driven_III}
F : \Omega \to \Omega, \quad \left( \begin{array}{c} \theta_{n+1} \\ x_{n+1} \end{array}\right) = F\left( \begin{array}{c} \theta_{n} \\ x_{n} \end{array}\right) := \left( \begin{array}{c} \theta_n + \vec\omega \\ g(\theta_n, x_n) \end{array}\right) .
\end{equation}
The space $\Omega$ is unknown, and its points are pairs of points from $\mathbb{T}^d$ and $\mathcal{X}$ respectively. This abstract formulation applies to all quasiperiodically driven systems. We further assume the dynamics in \eqref{eqn:def:quasi_driven} to be of the form
\begin{equation}\begin{split} \label{eqn:def:quasi_driven_II}
\theta_{n+1} &= \theta_n + \Delta t \vec\omega \bmod 2\pi \\
x_{n+1} &= g_{per}( \theta_n) + g_{chaos}(\theta_n, x_n)
\end{split}\end{equation}
Thus the task is to find (i) the \emph{quasiperiodicity dimension} $d$ and the rotation vector $\omega$; and the functions (ii) $g_{per}$ and (iii) $g_{chaos}$. The functions and spaces above will be assumed to be unknown. The only information about the system will be through a collection of $k$ observations / measurements, represented collectively as a map $Y : \Omega \to \real^k$. $Y$ is possibly unknown, and possibly a low-dimensional / partial observation of $\Omega$. It generates a sequence of $k$-dimensional data points $\SetDef{ y_n := Y(\theta_n, x_n) }{ n=0,1,2,\ldots }$, where $(\theta_n, x_n)$ is a trajectory of the dynamics in \eqref{eqn:def:quasi_driven_III} under $F$.

\subsection{Koopman operator and its spectrum} \label{sec:eigen}

The Koopman operator $U$ is essentially a time-shift operator. It operates on functions instead of points on the phase space. Given a function $\phi:\Omega \to \real$, $U\phi$ is the function defined as
\begin{equation} \label{eqn:def:Koop}
(U\phi)(z) := \phi\left( F z\right) , \quad \forall z\in\Omega ,
\end{equation}
%
%\[\begin{split} & U(\phi_1 + \phi_2)(z)  = (\phi_1 + \phi_2)(Fz) = \phi_1(Fz) + \phi_2(Fz) = (U\phi_1 + U\phi_2)(z) \end{split}\]
%
where $F$ is the underlying dynamical system \eqref{eqn:def:quasi_driven_III}. $\phi$ can be interpreted as a measurement or observation on the phase space $\Omega$, and $U\phi$ is the evolution/transformation of this measurement with the dynamics.

\paragraph{Koopman eigenfrequencies} Since $U$ is unitary; its spectrum must lie on the unit circle of the complex plane. The eigenvalues of $U$ correspond to the point spectrum, and any eigenfunction $\zeta$ has a corresponding eigenvalue of the form $e^{\iota \omega}$ for some $\omega\in \real$. Thus
\begin{equation} \label{eqn:Koop_eig}
(U^n \zeta)(z) \stackrel{\mbox{by def.}}{=} \zeta( F^n z) = e^{\iota \omega n} \zeta(z) , \quad \forall n\in\num .
\end{equation}
$\omega$ is called the Koopman eigenfrequency corresponding to $\zeta$. Thus the time-evolution of Koopman eigenfunctions is equivalent to multiplication by $e^{\iota \omega n}$ as a function of time $n$. 

$U$ always has the constant functions as eigenfunctions with eigenfrequency $0$. $U$ may or may not have other eigenfrequencies. The collection of eigenfunctions and (eigen)-frequencies have an algebraic structure to them. For any two  frequencies $\omega_1, \omega_2$, and integers $a,b$, $a\omega_1 + b\omega_2$ is also a frequency \cite{DasGiannakis_RKHS_2018}. If the system has at least one nonzero frequency, then it has all harmonics of that frequency and thus infinitely many frequencies. A collection of eigenfrequencies is said to be independent if no integer linear combination is an integer.
If the system has two independent frequencies, then all its frequencies are together dense on the real line. A collection of frequencies will be called a \emph{basis} or \emph{generating} set of eigenfrequencies if they are independent, and all frequencies of the system can be generated by taking integer linear combinations of frequencies from this set. There is no unique choice of a basis, but all bases will have the same cardinality $d$, called the quasiperiodicity dimension $d$. $d$ is a fixed finite number if $\Omega$ is a finite-dimensional manifold.

\paragraph{Koopman eigenfunctions and torus dynamics} Koopman eigenfunctions reveal quasiperiodic dynamics embedded in the system. Any $k$ Koopman eigenfunctions leads to a rotation on a $k$-dimensional torus:
\[\begin{tikzcd}
\Omega \arrow{d}[swap]{\zeta_1, \ldots, \zeta_k} \arrow{r}{F} & \Omega  \arrow{d}{\zeta_1, \ldots, \zeta_k} \\
\mathbb{T}^k \arrow{r}{R_{\vec\omega}} & \mathbb{T}^k
\end{tikzcd}\]
Here, $R_{\vec\omega}$ is a rotation by the vector $\vec\omega$ of frequencies. If these eigenfrequencies are independent, then this map will be surjective. Taking $d=k$ implies that the dynamics has an embedded / factor torus rotation of the same dimension as the quasiperiodicity dimension.

This completes our examination of the quasiperiodic structure of the dynamics \eqref{eqn:def:quasi_driven_III}. We next discuss some techniques from Functional Analysis for reconstructing the quasiperiodic component and its complement.

%-_-_-_-_-_-_-_-_-_-_-_-_-_-_-_-_-_-_-_-_-_-_-_-_-_-_-_-_-_-_-_-_-_-_-_-_-_-_-_-_-_-_-_-_-_-_-_-_-_-_-_-_-_-_-_-_-_-_-_-_-_-_-_-_-_-_-_-_-_-_-_-_-_-_-_-
\subsection{Kernels and integral operators} \label{sec:kernel}

A kernel is a function $ k : M\times M\to \real $ on some space $M$. The quantity $k(x,y)$ is measure of similarity, closeness or distance between two points $x,y\in M$.  Kernel based methods have been used very effectively to obtain information such as statistical manifolds  \cite{DasDimitEnik2020}, geometric information \cite{BerrySauer2017}, and dynamical information such as tracer flows \cite{GiannakisDas_tracers_2019}, stable/unstable foliations \cite{BerryEtAl2013}, Koopman spectrum  \cite{DasGiannakis_delay_2019, DasGiannakis_RKHS_2018}. The techniques in this paper are based on \cite{DasGiannakis_RKHS_2018}. We shall use the \emph{Gaussian kernel}
\[ k_\epsilon(x,y) := \exp\left( -\frac{1}{\epsilon} d(x,y)^2 \right), \]
where $\epsilon$ is called the \emph{bandwidth parameter}, and $d(\cdot, \cdot)$ is some notion of metric or distance on the space.

\paragraph{Delay-coordinates} The data sequence $y_n$ is obtained through an observation $Y$. However $Y$ may not be a one-to-one map and its values may not correspond to unique states in $\Omega$. We convert $Y$ into an embedding using the method of delay coordinates \cite{SauerEtAl91}, by incorporating $Q$ delays to get the map $Y^{(Q)} : \Omega \to \real^{k(Q+1)}$ :
\[Y^{(Q)}(\omega) := \left( Y(\omega), Y( F^1 \omega), \ldots, Y(F^Q\omega) \right) . \]
Thus the delay coordinated version of each point $y_n$ is
\[ y_n \leftrightarrow y_N^{(Q)} := \left( y_n y_{n+1}, \ldots, y_{n+Q} \right) . \]
We next use the Gaussian shape function to implicitly obtain a kernel $k_\epsilon : \Omega \times \Omega \to \real$ as follows 
\begin{equation} \label{eqn:def:Gauss_ker} 
k_\epsilon( z, z') := \exp\left( -\frac{1}{\epsilon} \norm{ Y^{(Q)}(z) - Y^{(Q)} (z') }^2  \right).
\end{equation}

Even if the two states $z, z'$ are unknown, the left-hand side in \eqref{eqn:def:Gauss_ker} can be computed since the right-hand side only uses the observation map $Y$. We next modify $k_\epsilon$ by a process called \emph{bistochastic normalization} \cite{DGJ_compactV_2018} to get a kernel $p_\epsilon$ which is symmetric, Markovian, and more adapted to the non-uniform distribution of the data.
\[\begin{split}
& \deg_R(z) := \int k_\epsilon( z, z') d\mu(z'), \, \deg_L(z) := \int \frac{ k_\epsilon( z, z') }{ \deg_R(z) } d\mu(z') . \\
&  \tilde{k}_\epsilon(z, z'') := \frac{k_\epsilon(z, z'')}{ \deg_R(z) \deg_L(z'')^{1/2} } \\
& p_\epsilon(z, z') := \int \tilde{k}_\epsilon(z, z'') \tilde{k}_\epsilon(z'', z') d\mu(z'') .
\end{split}\]

\paragraph{Kernel integral operator} Associated to the kernel $p_\epsilon : \Omega \times \Omega \to \real$ is the integral operator $P_\epsilon$, which operates on $L^2(\mu)$ functions $\phi:\Omega\to \real$ as
\[(P_\epsilon  \phi)(z) := \int_{\Omega} p_\epsilon(z, z') \phi(z') d\mu(z').  \]
$ P_\epsilon$ is a  \emph{compact}, symmetric operator on $L^2(\mu)$ \cite[e,g,][]{DGJ_compactV_2018}. Moreover $ P_\epsilon$ has a complete basis of eigenfunctions
\[  P_\epsilon \phi_j = \lambda_j \phi_j , \quad j=0,1,2,\ldots , \]
where the indexing is done so that the $\lambda_j$s are in decreasing order. Due to the normalizations carried out, we have $\phi_1 \equiv 1_\Omega$, the constant function equal to $1$ everywhere. Moreover,  the eigenvalues satisfy$1 = \lambda_1 \geq \lambda_2 \geq \lambda_2 \geq \ldots > 0$.  Also importantly, the $\phi_j$ are an orthonormal basis, i.e.,
\[ \langle \phi_i, \phi_j \rangle_{L^2(\mu)} := \int \phi_i^*(x) \phi_j(x) d\mu(x) = \delta_{i,j} . \]
All these properties of the $\lambda_j$ and $\phi_j$ are useful for \emph{kernel-based learning}, in which we recreate or extrapolate unknown functions from some samples using these $\phi_j$s as a basis.

\paragraph{Kernel based learning} One of the main advantages of kernel-based approaches is that while the $\phi_j$ can be approximated to any degree of accuracy by solving an eigenvalue equation of a data-driven matrix; they can be easily extended from vectors to a continuous function over the entire data space $\real^{k(|Q|+1)}$ and thus on $\Omega$. The $\phi_j$s also happen to be left singular vectors of the asymmetric operator $\tilde{K}_\epsilon$, with $\sqrt{\lambda_j}$ and $\gamma_j$ being the associated singular values and right singular vectors. Thus for an arbitrary point $z\in \real^{k(|Q|+1)}$, we have
\begin{equation} \label{eqn:dskn30}
\begin{multlined}
\phi_j(z) = \lambda_j^{-1}\int p_\epsilon(z , z') \phi_j(z') d\mu(z') \\= \lambda_j^{-1}\int \tilde{k}_\epsilon(z , z') \gamma_j(z') d\mu(z') .
\end{multlined}
\end{equation}
These $\phi_j$ for the basis for \emph{learning} any function $f : \real^{k(Q+1)} \to \real^d$. It is done by first computing the components of $f$ along the first $L$  $\phi_j$
\[ f_l := \langle \phi_l, f \rangle_{L^2(\mu)} = \int_{\real^{k(Q+1)}} \phi_l^*(z) f(z) d\mu(z) ,  \]
and then taking the sum/integral
\[  f(y) \approx \sum_{l=1}^{L} f_l \phi_l(y) = \int_{\real^{k(Q+1)}} p_\epsilon ( y, z) \sum_{l=1}^{L} \frac{ f_l }{ \lambda_l } \phi_l(z) d\mu(z) .  \]
The parameter $L$ is called the \emph{spectral truncation} parameter; it is the size of the hypothesis space. 

We next describe the use of the fast-Fourier transform to derive the Koopman eigenfrequencies but on the $L$ eigenfunctions instead of the raw data. 

%-_-_-_-_-_-_-_-_-_-_-_-_-_-_-_-_-_-_-_-_-_-_-_-_-_-_-_-_-_-_-_-_-_-_-_-_-_-_-_-_-_-_-_-_-_-_-_-_-_-_-_-_-_-_-_-_-_-_-_-_-_-_-_-_-_-_-_-_-_-_-_-_-_-_-_-
\section{The data-driven procedure} \label{sec:data_driven_algorithms}

In the data-driven approach, all of the entities described in Section~\ref{sec:kernel} have a \emph{data-driven analog}. We begin with the invariant measure $\mu$ itself. It will be replaced by $\mu_N = \frac{1}{N} \sum_{n=1}^{N} \delta_{y_n^{(Q)}}$, the average of the Dirac-delta measure on the data points $y_n^{(Q)}$. These are called \emph{sampling/empirical} measures. Their integrals with respect to these sampling measures are given by
\[  \int_{\real^{k(Q+1)}} \phi d \mu_N = \frac{1}{N} \sum_{n=1}^{N} \phi\left(  \delta_{y_n^{(Q)}} \right) , \]
for every continuous test function $\phi : \real^{k(Q+1)} \to \real$. The kernel integral operators $K$, $\tilde{K}$ and $P$ will be approximated as $N\times N$ matrices $\Matrix{K}$, $\Matrix{\tilde K}$ and $\Matrix{P}$ :
\begin{equation} \label{eqn:def:K_data}
\Matrix{K}_{i,j} := k_\epsilon\left( y_i^{(Q)}, y_j^{(Q)} \right), \; \vec{d} := \frac{1}{N} \Matrix{K} \vec{1}_N, \; .
\end{equation}
Next compute 
\begin{equation} \label{eqn:def:vecq}
\vec{q} := \frac{1}{N} \Matrix{K} D^{-1}, \quad D:= \diag( \vec{d} ) .
\end{equation}
$\vec{d}$ and $\vec{q}$ are called \emph{left and right degree vectors} respectively. Finally compute 
\begin{equation} \label{eqn:def:tildeK_data}
 \Matrix{\tilde{K}} := D^{-1} \Matrix{K} Q^{-1/2}, \quad Q := \diag( \vec{q} ) .
\end{equation}
The matrix $\Matrix{P} := \Matrix{\tilde{K}} \Matrix{\tilde{K}}^*$ is not computed explicitly. 
See \cite[Algorithm 1]{DGJ_compactV_2018} for more details.  Compute the top $L$ singular values $1=\sigma_1 \geq \ldots \geq \sigma_L$ of $\Matrix{\tilde{K}}$ and the corresponding left eigenvectors $ \vec{\phi}_{1}, \ldots,  \vec{\phi}_{L}$ and right singular vectors $\vec{\gamma}_{1}, \ldots, \vec{\gamma}_{L}$. Set $\lambda_i := \sigma_i^2$, for $i=1,\ldots, L$.
The set of vectors $\SetDef{ \vec{\phi}_{l} }{ l = 1, \ldots, L }$ and $\SetDef{ \vec{\gamma}_{l} }{ l = 1, \ldots, L }$ are both orthonormal systems for $\cmplx^N$. The eigenvectors have continuous extensions to any $y\in \real^{k(Q+1)}$ as
\begin{equation} \label{eqn:oss_data}
\begin{split}
& \bar{\phi}_l (y) := \frac{1}{N \lambda_l} \vec{ k }_\epsilon (y)^\top Q^{-1/2} \vec\gamma_l \\
& \vec{ k }_\epsilon (y) := \left( k_\epsilon \left( y, y_1^{(Q)} \right)  , \ldots, k_\epsilon \left( y, y_N^{(Q)} \right) \right) .
\end{split}
\end{equation}
The function $ \bar{\phi}_l$ is continuous as the vector $ \vec{k}_\epsilon (y) $ is a continuous function of $y$. If $y$ in the above equation is substituted by one of the data-points $y_n^{(Q)}$, then by design $\bar{\phi}_l \left( y_n^{(Q)} \right) = \vec{\phi}_{l,n}$.  Thus $ \bar{\phi}_l$ is indeed a continuous extension of the vector $\vec{\phi}_l$. This feature of extendability and easy evaluation at arbitrary points is one of the most powerful tools of kernel-based methods. %In our next algorithm, we show a different application of these eigenfunctions, for discovering true Koopman eigenfrquencies.

\paragraph{RKHS based spectral filtering} The following procedure was described in \cite[Algorithm 1]{DasGiannakis_RKHS_2018}, and accepts as parameters  $\epsilon_1, \epsilon_2>0$ and integer $L_0>1$. Let $\Fourier_N$ denote the discrete Fourier transform on vectors of length $N$. Let $\Matrix{\Phi}$ be the $N\times L$ matrix whose $l$-th column is $\vec{\phi}_{l}$. Set $\Lambda:= \diag\left(  \lambda_1, \ldots,\lambda_L \right)$ and compute
\[ \Matrix{\hat\Phi} := \Fourier \Matrix{\Phi} \quad \Matrix{H} := \Matrix{\hat\Phi} \Lambda^{-0.5}. \]
Next, compute the RKHS-norms as
\[ \Matrix{W}_{n,1} := \abs{ \Matrix{H}_{n,1} },\] \[\Matrix{W}_{n,l+1} := \Matrix{W}_{n,l} + \abs{ \Matrix{H}_{n,l+1} }, l=1,\ldots, L-1 .  \]
Now let $J = 1,\ldots, N$. Discard all the $j\in J$ for which $\Matrix{W}_{j,L_0}<\epsilon_1$. Of the remaining $j\in J$, discard those $j$ for which $\ln \Matrix{W}_{j,L} - \ln \Matrix{W}_{j,L_0} > \epsilon_2$.  Compute $\omega_j = \frac{2\pi j}{N \Delta t}$ for all of the remaining $j\in J$. These frequencies $0=\omega_1<\omega_2<\ldots<\omega_m$ can be interpreted to be true Koopman eigenfrequencies with substantial presence in the original data.

\paragraph{Periodic and chaotic components} The identified frequencies $0=\omega_1<\omega_2<\ldots<\omega_m$ are by no means exhaustive, they are only a finite subset of a usually infinite set of Koopman eigenfrequencies. However, they represent those (true) frequencies that have a significant presence in the data. The threshold $\epsilon_1$ is meant to be a numerical implementation of frequencies being significant. We next
construct the periodic component $g_{per} : \real \to \real^k$ as
\begin{equation} \label{eqn:def:Acoeff}
g_{per}(t) := \Re \sum_{j=1}^{m} \Matrix{A}_{j,:} e^{\iota \omega_j t} .
\end{equation}
The $m\times k$ matrix $\Matrix{A}$ is the least-squares solution to 
\[ \Re \Matrix{F} \Matrix{A} = \Matrix{Y}, \]
where $Y$ is the data-matrix and $\Matrix{F}$ is $N\times m$ matrix
\[ \Matrix{F}_{n,j} := (2 - \delta_{j,1}) e^{\iota n \omega_j}  \quad 1\leq n\leq N, \, 1\leq j\leq m .  \]
Next, we construct the chaotic component 
\begin{equation} \label{eqn:def:Ecoeff}
    g_{chaos} : \real^{k(Q+1)} \to \real^k, \; g_{chaos} := \sum_{l=1}^{L} \Matrix{E}_{l,:} \bar{\phi}_l ,
\end{equation}
where $\Matrix{Y_{non}} :=  \Matrix{Y} - \Re \Matrix{F} \Matrix{A}$ and $\Matrix{E} := \Matrix{\Phi}^* \Matrix{Y_{non}}$.

\paragraph{The reconstruction} We avoid the task of identifying a set of generating frequencies by directly using the selected frequencies in the approximation
\begin{equation} \label{eqn:sdkjn39}
g_{per} \left( \vec\theta + n\vec{ \tilde{\omega} } \right) \approx \sum_{j=1}^{m} a_{j} \exp\left( \iota n \omega_j  \right) ., \quad n=0,1,2,\ldots . 
\end{equation}
Using this simplification in~\eqref{eqn:sdkjn39}, and the formulas in \eqref{eqn:def:Ecoeff} and \eqref{eqn:def:Acoeff}, we create the following data-driven model of the dynamics :
\begin{equation} \label{eqn:reconstruct} 
\begin{split}
y^0_{n+1} &:=  g_{per}(n\Delta t) + g_{chaos}\left( y^0_n , \ldots, y^Q_n   \right) \\
y^1_{n+1} &:= y^0_n \\
\vdots &= \vdots \\
y^Q_{n+1} &:= y^{Q-1}_n 
\end{split} .
\end{equation}

Here each $y^q_n \in \real^k$, thus making the state vector $y_n = \left( y^0_n , \ldots, y^Q_n \right)$ a vector in $\real^{k(|Q|+1)}$. We have thus created a standalone dynamical system $\real^{k(|Q|+1)}$ which is conjugate to the latent dynamics.

\section{Case Study: Signalized Intersection Corridor} \label{sec:Alfaya}
\subsection{Data Description}

The case study analyzes queue length measurements from nine adaptive traffic signals located on the Alafaya Trail (SR-434) in East Orlando, FL. The obtained data includes the details of each movement with the time, duration, queue length, and waiting time. It provides information on eight movements: north left (NL), north through (NT), south left (SL), south through (ST), east left (EL), east through (ET), west left (WL), and west through (WT). In this study, we focus on the queue length formation of northbound through movements. The raw data was processed and calibrated by Rahman, and et al. \cite{rahman2021real} and was resampled at regular intervals of $\Delta t = 2 \text{minutes}$. This work uses the processed data from \cite{rahman2021real}.

\begin{figure}[!htbp]\center
\includegraphics[width=1.0\linewidth]{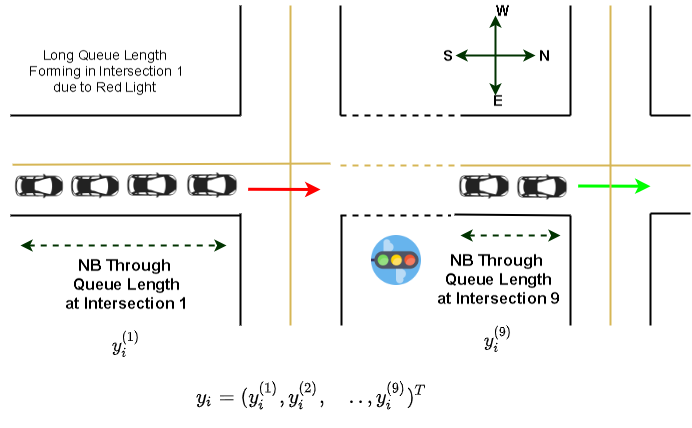}
\vspace{-2mm}
\caption{Northbound queue length formation at the corridor at $i^{th}$ time}
\label{fig:dgh3}
\end{figure}

\subsection{Problem Formulation}
We formulate the signalized intersection corridor as a quasiperiodically driven dynamical system. The underlying dynamics of the system are high dimensional, and its governing equations are unknown. We hypothesize that the signalized intersection corridor system obeys a dynamics of the form \eqref{eqn:def:quasi_driven_III}, and the observed queue lengths are generated through some measurement function $Y$, as described in Section~\ref{dynamical_theory}. Note that the measurement $Y$ is not necessarily one-to-one. In particular, it may not be possible to connect the $y_i$ with a dynamical rule of the form $y_{i+1} = \tilde{F}(y_i)$. Rather $y_i$ should be interpreted as a partial observation of the true state in $\Omega$, in the $i$-th time frame. In Figure \ref{fig:dgh3} we show the visual representation of $y_i$. We use only this data to obtain a parameter-free  reconstruction  of  the dynamical system.

\subsection{Computation}
\begin{figure}[!htbp]\center
\includegraphics[width=.49\linewidth]{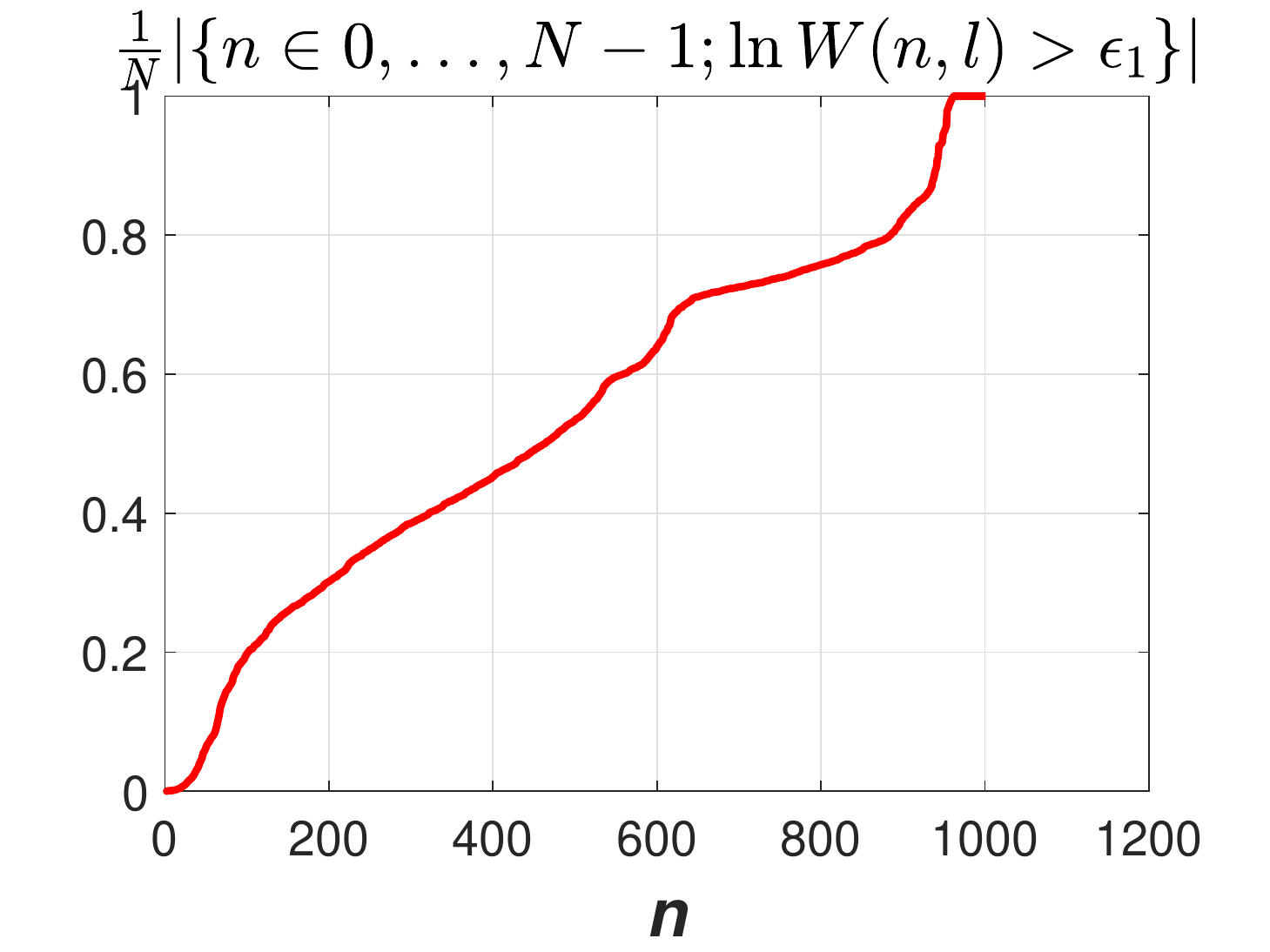}
\includegraphics[width=.49\linewidth]{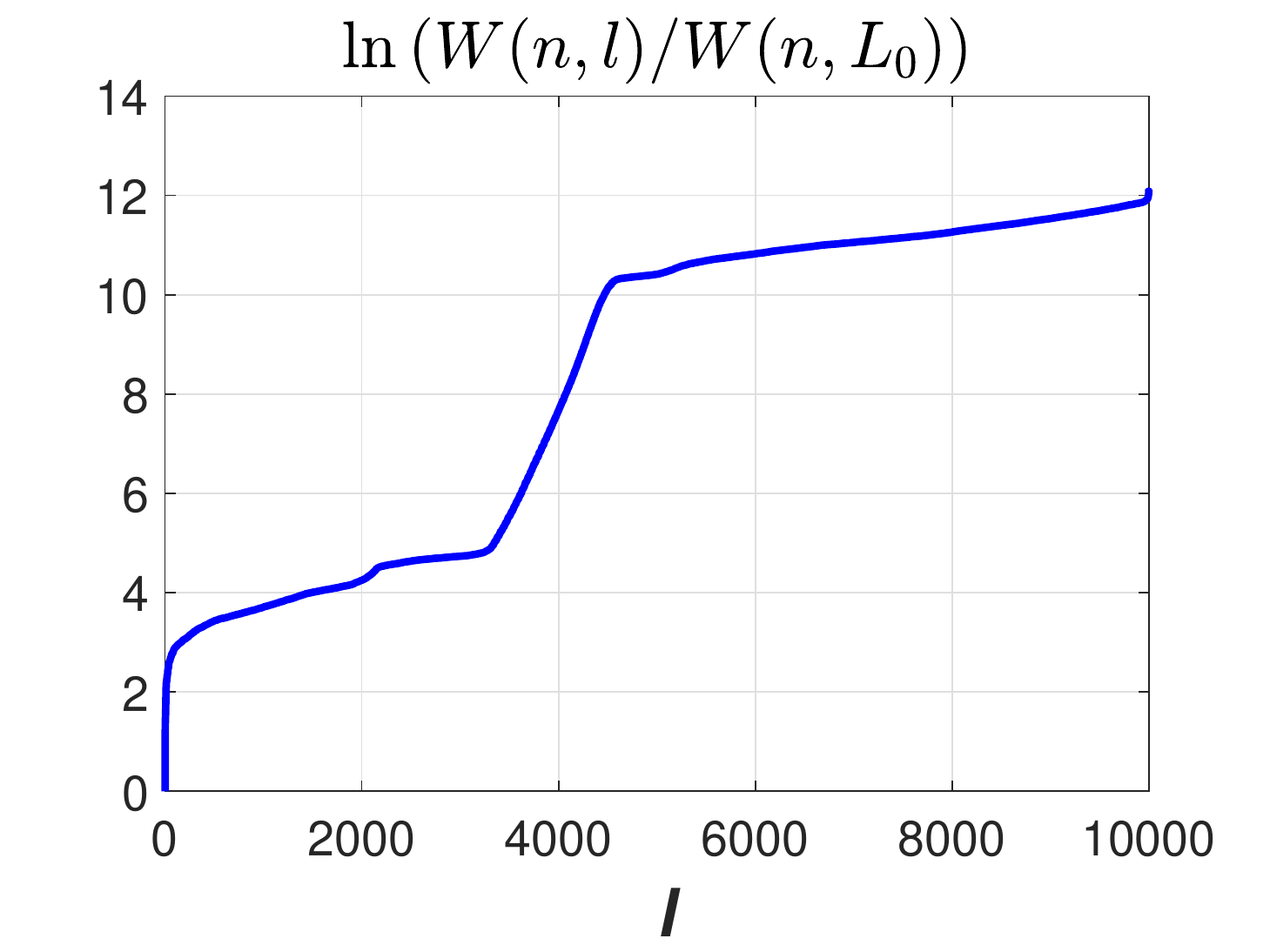}
\vspace{-2mm}
\caption{ Choice of the thresholds $L_0, \epsilon_2$ for the traffic intersection data. $L_0$ was chosen to be $100$, as it corresponds to the approximate location of an abrupt jump in the curve in the left panel. Given this choice of $L_0$, we plotted the quantity $\ln \left( W(n,l) / W(n,L_0) \right)$ as a function of $n \in 0,\ldots, N/2$. We choose $\epsilon_2$ to be $2.5$, which is approximately the first of the two points of inflection of the curve in the red panel, as indicated by the horizontal red and blue lines. }
\label{fig:eps_12}
\end{figure}

We arrange the data in the form of a matrix $\Matrix{Y}$ with $k=9$ columns. Each column corresponds to the traffic queue length as a function of time, at one among $9$ intersections along Alfaya Trail. We used a bandwidth parameter of $\epsilon = 0.1$ for the Gaussian kernel. We compute a total of $L=1001$ eigenfunctions for the bistochastic kernel. We set $\epsilon_1 = 0.1$ and then choose the parameters $L_0=5$ and $\epsilon_2 = 2.5$ using the heuristic approach shown in Figure~\ref{fig:eps_12}. 
%The left panel shows the fraction $\frac{1}{N} \SetDef{  n\in 1,\ldots, N }{  W[n,l]>\epsilon_1 } $ as a function of  $l$, as $l$ varies from $1$ to $L=1001$. A good choice of $L_0$ would be that value of $l$ at which the graph becomes steep. There is no unique such $L_0$ and we choose $L_0=100$. The second panel is an aid to choose the other threshold $\epsilon_2$. We collected the values  $\SetDef{ \ln \Matrix{W}_{n,L} - \ln \Matrix{W}_{n,L_0} }{ n\in 1,\ldots,N }$ and the graph shows the distribution of these values. A good choice of $\epsilon_2$ would be a point where the graph becomes steep. 
Lower the value of $\epsilon_2$, more frequencies get filtered out and higher the probability of the identification being correct. The two thresholds $\epsilon_1$ and $\epsilon_2$ are based on the asymptotic behavior in two different directions \cite[Theorem 1,4]{DasGiannakis_RKHS_2018}. Combined, they provide a surer guarantee of identification of true eigenfrequencies and the discarding of \emph{spurious} eigenvalues or \emph{pseudo}-spectrum.

\subsection{Results}

\subsubsection{Generating Frequencies}
\begin{figure}[!htbp]\center
\includegraphics[width=1.0\linewidth]{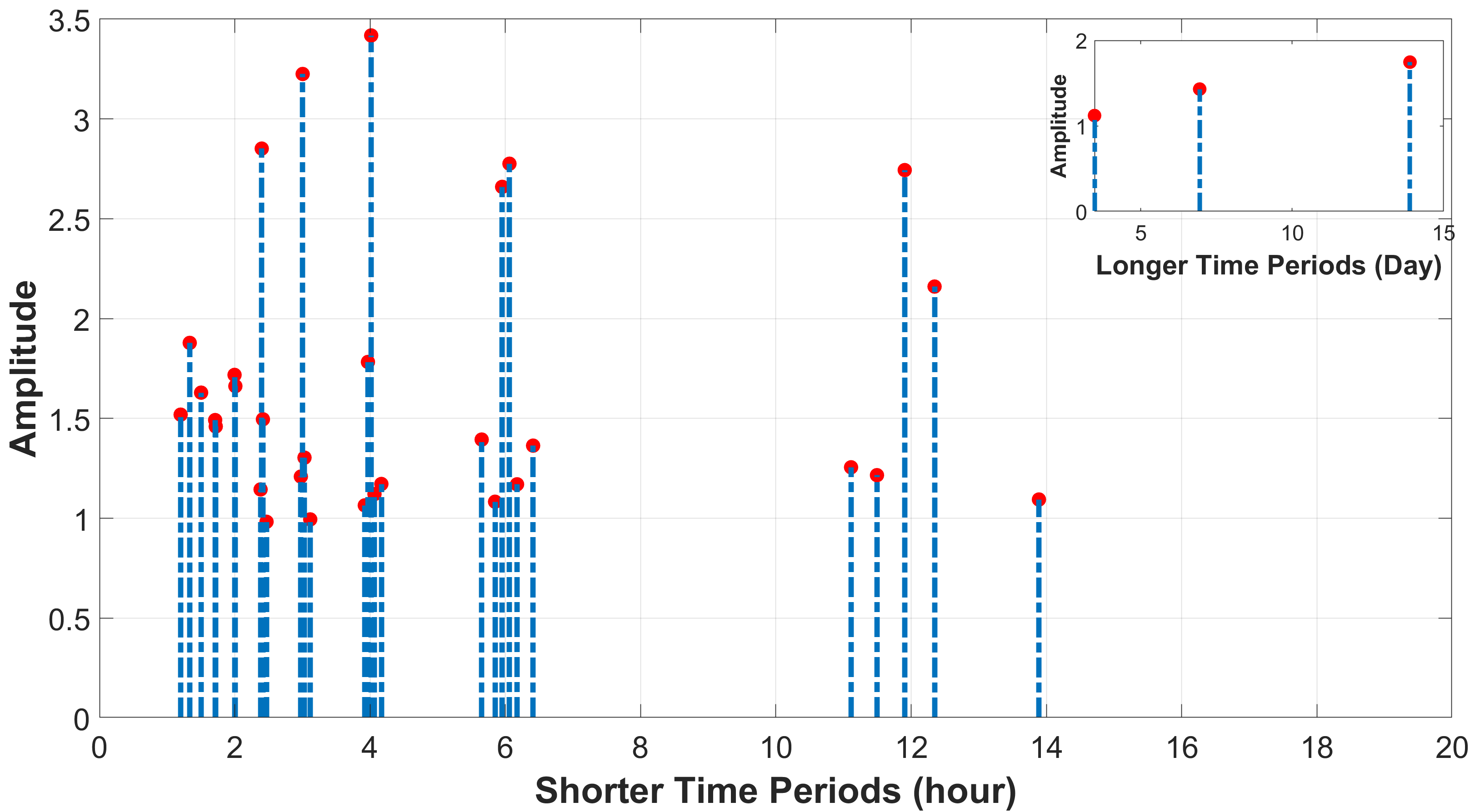}
\vspace{-2mm}
\caption{ Selected frequencies and time-periods.}
\label{fig:freq_sel}
\end{figure}

This work analyzed $20,000$ snapshots of queue lengths for the nine intersections to extract Koopman frequencies, i.e., generating frequencies of the quasiperiodic driving source of the dynamics. Figure~\ref{fig:freq_sel} exhibits the spectrum of Koopman eigenfrequencies. The $x$ axis denotes periods corresponding to the frequencies. We separate longer periods from the shorter ones and present them in two different panels for convenience of presentation.  For both the panels, the y-axis shows $\Matrix{W}_{j, L_0}$ which we have denoted as amplitude for the selected period indices $j$. Figure~\ref{fig:freq_sel} shows that dominant periods are clustered around 1 hour, 2 hours, 3 hours, 6 hours, 12 hours, 14 hours, 3.5 days, 7 days, and 14 days. The identified periods correspond to the natural periods of the system. These periods are consistent with the results in \cite{wang2014multiscale}. The authors decomposed speed measurements from an intersection corridor via multiscale multifractal analysis (MMA). They reported that the dominant periodicities on weekdays are 7 days, 24, 12, 8, 6, and 3 hours while on weekends are 12 and 24 hours. In the present work, we did not differentiate weekend data from weekday data. However, we took a different set of observables (i.e., northbound through queue length data) and a completely different intersection corridor system, our identified quasiperiodic frequencies matched with that of \cite{wang2014multiscale}. This finding corroborates that the proposed technique can successfully identify generating frequencies of the quasiperiodic driving force of the intersection dynamical systems. 

\subsubsection{Reconstruction and Prediction}

\begin{figure}\center
\includegraphics[width=1\linewidth]{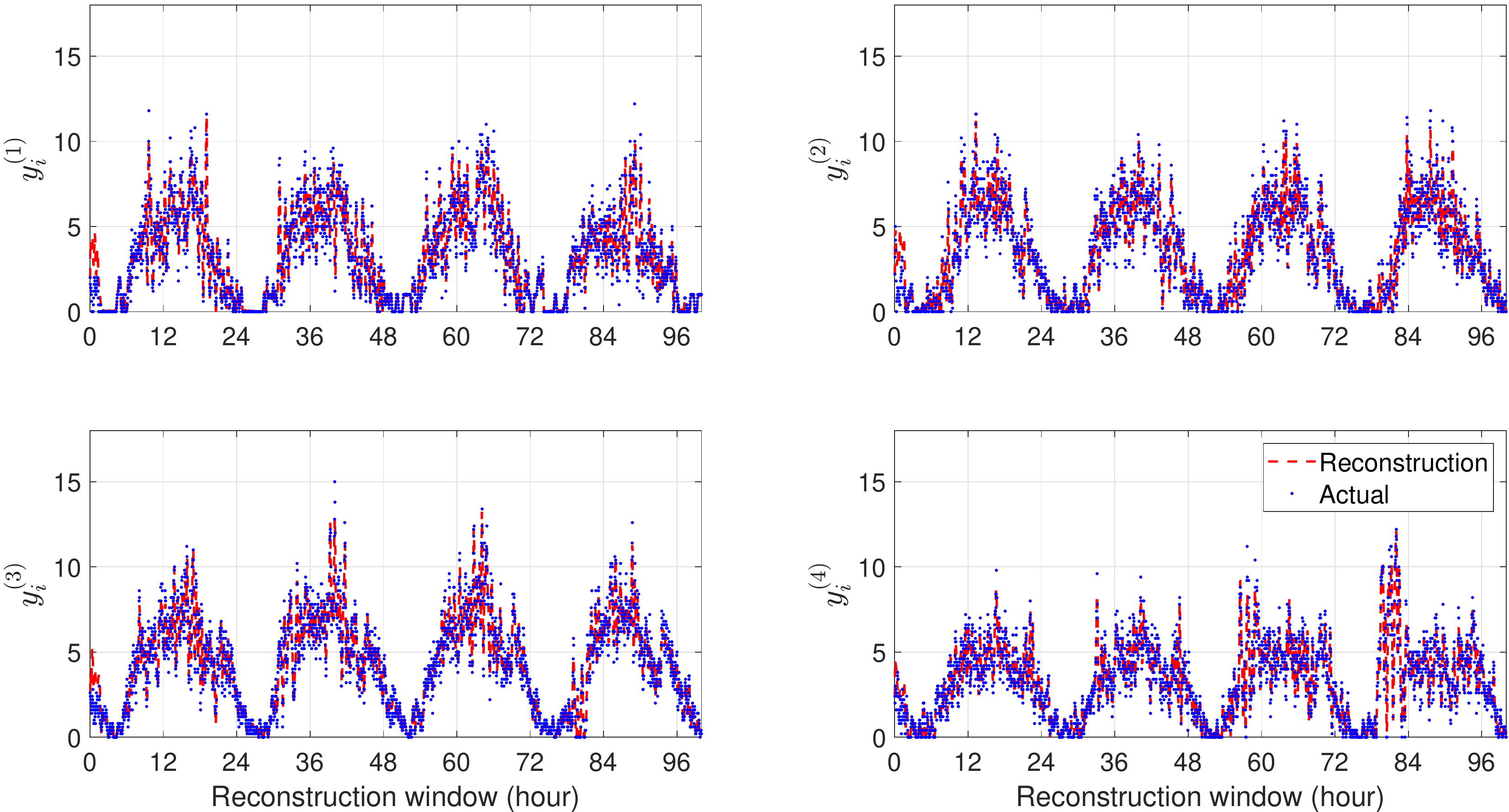}
\vspace{-2mm}
\caption{Dynamical reconstruction of NB queue lengths, shown for the first $4$ intersections.}
\label{fig:reconstruct}
\end{figure}

\begin{figure}\center
\includegraphics[width=1\linewidth]{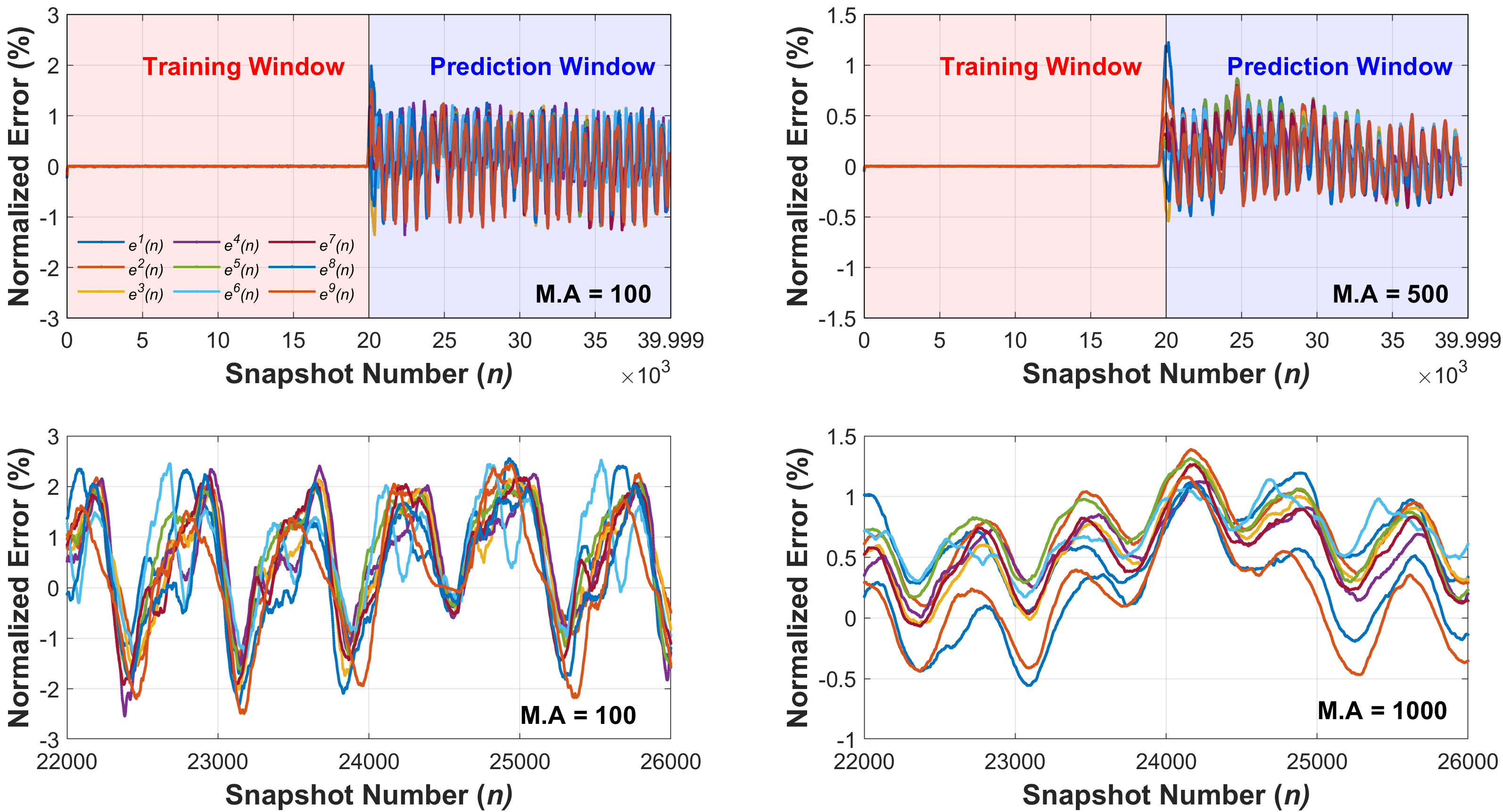}
\vspace{-2mm}
\caption{ Dynamical reconstruction and relative prediction error \eqref{eqn:def:err_rel} of NB queue length data. Normalized error for intersections in $n^{th}$ snapshot are indicated as $e^1(n), \ldots, e^9(n)$. The top two plots shows a sudden and qualitative difference in the error beyond the training period of $n=20000$. The bottom two plots focus on the prediction performance beyond the training period. The system \eqref{eqn:reconstruct} initialized with the state $y^{(Q)}_{t=2200}$, and iterated over the time window $n = [22000, 26000]$. An increased size of the moving average window reduces the effect of random outliers or events and diminishes the prediction error. Number of snapshots used to calculate moving averages in each cases is denoted as M.A and shown in each plot.}
\label{fig:pred}
\end{figure}

We reconstruct the original data from the decomposed parts, which is shown in Figure~\ref{fig:reconstruct}. The red curve is the output obtained from the reconstruction at each intersection. The curves closely following each other, including the moments when fluctuations occur. Although reconstructed dynamical models differ from the true system, this difference is inevitable in a learning problem. However, the reconstruction of $g_{per}$ is bounded, which guarantees that the dynamics under \eqref{eqn:def:quasi_driven_II} would remain bounded, and the deviation of the trajectories also remain bounded. In Figure~\ref{fig:pred} we illustrate the error in reconstruction and prediction by computing the normalized relative error 
\begin{equation} \label{eqn:def:err_rel}
\text{error}_{rel}^{(i)}(n) := \frac{ \abs{ y_n^{(i)} - \hat{y}_n^{(i)} } }{ \max_n \abs{ y_n^{(i)} } }, i = 1,\ldots 9.
\end{equation}
Here $\hat{y}_n^{(i)}$ denotes the output of the reconstructed system \eqref{eqn:reconstruct}.        

\section{Conclusion}
This work developed a data-driven framework for modeling quasi-periodic dynamical systems using Koopman theoretic approach. The proposed approach can handle dynamics with strong nonlinear and chaotic components. Thus, it is applicable across domains. We performed a case study using queue length data on a corridor of nine signalized intersections. The proposed approach accurately identified the generating frequencies and the results for reconstruction and prediction are encouraging. The long-term prediction error remained bounded without exogenous inputs, unlike recurrent neural network-based methods such as long short-term memory (LSTM). Moreover, in comparison to deep NNs, the proposed technique is not a black-box approach. All these advantages make it a promising candidate for future research.   

\section*{Acknowledgements}
Authors would like to thank Dr. Samiul Hasan and his research group for providing access to the intersection dataset. 
\printbibliography
\end{document}